\documentclass[1p]{elsarticle}
\usepackage{amsfonts, amsmath, amssymb, amsthm}
\usepackage[utf8]{inputenc}
\usepackage{lineno, blindtext}
\usepackage{array, makecell}
\usepackage{tikz, color, xcolor, soul}
\usepackage{pgfplots}
\usepackage{bbold}
\usepackage{enumitem}
\newtheorem{theorem}{Theorem}

\newtheorem{lemma}[theorem]{Lemma}

\usepackage{ulem}

\theoremstyle{definition}
\newtheorem{definition}{Definition}
\newtheorem{remark}{Remark}
\newtheorem{construction}{Construction}

\newcommand{\F}{\mathbb{F}}
\newcommand{\Fq}{\mathbb{F}_q}

\newcommand{\cC}{{\mathcal C}}
\newcommand{\cF}{{\mathcal F}}

\renewcommand{\u}{\boldsymbol{u}}
\newcommand{\e}{\boldsymbol{e}}
\renewcommand{\c}{\boldsymbol{c}}
\renewcommand{\v}{\boldsymbol{v}}

\newcommand{\w}{\boldsymbol{w}}

\newcommand{\x}{\boldsymbol{x}}
\newcommand{\y}{\boldsymbol{y}}

\newcommand{\0}{\boldsymbol{0}}

\newcommand{\aalpha}{\boldsymbol{\alpha}}
\newcommand{\bbeta}{\boldsymbol{\beta}}
\newcommand{\ind}{\mathbb{1}}

\DeclareMathOperator{\spn}{span}

\newif\ifcomment
\commenttrue

\newcommand{\new}[1]{\textcolor{black}{#1}}

\begin{document}

\begin{frontmatter}

\title{Almost Affinely Disjoint Subspaces}

\author[TUM]{Hedongliang Liu}
\author[TUM,Sk]{Nikita Polianskii}
\author[Sk]{Ilya Vorobyev}
\author[TUM]{Antonia Wachter-Zeh}

\address[TUM]{Technical University of Munich\\
                    Institute for Communications Engineering\\
                    80333 Munich, Germany\\
                    Email: \{lia.liu, nikita.polianskii, antonia.wachter-zeh\}@tum.de}
\address[Sk]{Skolkovo Institute of Science and Technology \\ Center for Computational and Data-Intensive Science and Engineering \\ 121205 Moscow, Russia }
\begin{abstract}
In this work, we introduce a natural notion concerning \new{finite vector} spaces.  A family of $k$-dimensional subspaces of $\F_q^n$\new{, which forms a partial spread,} is called almost affinely disjoint if any $(k+1)$-dimensional subspace containing a subspace from the family non-trivially intersects with only a few subspaces from the family. The central question discussed in the paper is the polynomial growth (in $q$) of the maximal cardinality of these families given the parameters $k$ and $n$. For the cases $k=1$ and $k=2$, optimal families are constructed. For other settings, we find lower and upper bounds on the polynomial growth.
Additionally, some connections with problems in coding theory are shown.
\end{abstract}


\end{frontmatter}

\section{Introduction}
Let $q$ be a prime power and let $\F_q$ denote the finite field with $q$ elements. By $\F_q^n$ we \new{denote} the standard vector space over $\F_q$, whose elements are $n$-tuples $\x=(x_1,\ldots,x_n)$ with $x_i\in\F_q$. A \textit{partial $k$-spread} in $\F_q^n$ is a collection of $k$-dimensional subspaces with pairwise trivial intersection.  For a vector $\u\in\F_q^n$ and a subspace $S\subseteq \F_q^n$, we define the affine subspace $\u+S:=\{\u+\v: \v\in S\}$. \new{Denote by $[m]$ the set of integers $\{1,\dots,m\}$.}

In this paper, we discuss two families of $k$-dimensional subspaces in $\F_q^n$ \new{satisfying} the following definitions.
\begin{definition}[Almost affinely disjoint subspace family]\label{def::super nice def}
  Given positive integers $k$ and $n$ such that $2k<n$, let $\cF$ be a family of $k$-dimensional linear subspaces in $\F_q^{n}$. This family is said to be \textit{$L$-almost affinely disjoint} (or, briefly, \textit{$[n,k,L]_q$-AAD}) if  the two properties hold:
  \begin{enumerate}
  \item The family $\cF$ is a partial $k$-spread of $\F_q^n$.
  \item For any $S\in \cF$ and $\u\in\F_q^n\setminus S$, the affine subspace $\u+S$ 
    intersects at most $L$ subspaces from the family $\cF$. 
  \end{enumerate}
\end{definition}
\new{Note that $n\geq 2k$ is required to construct a partial spread of size larger than $1$. We also exclude $n=2k$, since if the family $\cF$ is a partial $k$-spread of $\F_q^n$, the affine subspace $\u+S$, $S\in \cF$ and $\u\in\F_q^n\setminus S$, intersects every subspace from $\cF$ except $S$.}

\begin{definition}[Almost sparse subspace family]\label{def::sparse family}
	Given positive integers $k$ and $n$ such that $2k<n$, let $\cF$ be a family of $k$-dimensional linear subspaces in $\F_q^{n}$. This family is said to be \textit{$L$-almost sparse} (or, shortly, \textit{$[n,k,L]_q$-AS})  if  the two properties hold:
	\begin{enumerate}
	\item The family $\cF$ is a partial $k$-spread of $\F_q^n$.
		\item Any $(k+1)$-dimensional subspace in $\F_q^n$ intersects non-trivially at most $L$ subspaces from the family $\cF$.
	\end{enumerate}
\end{definition}
\begin{remark}
We also remark that the spread property does not affect the asymptotic analysis when $q\to\infty$ (see the proof of Theorem~\ref{th::upper bound AAD family}).
\end{remark}
Given $k$, $n$, $L$ and $q$ we denote the \textit{maximal size} of $[n,k,L]_q$-AAD and $[n,k,L]_q$-AS families by $m^{AAD}_q(n,k,L)$ and $m^{AS}_q(n,k,L)$.
Define the \textit{polynomial growth} of the maximal size of AAD and AS families by
\begin{align*}\label{eq::polynomial growth}
p^{AAD}(n,k,L)&:=\limsup_{q\to\infty} \log_q(m^{AAD}_q(n,k,L)),\\
p^{AS}(n,k,L)&:=\limsup_{q\to\infty} \log_q(m^{AS}_q(n,k,L)).
\end{align*}
\subsection{Related work}
The concept of almost sparse families is closely related to so-called \textit{subspace designs} introduced by Guruswami and Xing in~\cite{guruswami2013list} and further developed by Guruswami, Kopparty, Xing and Yuan in~\cite{guruswami2016explicit,guruswami2018subspace}. \new{A collection $\mathcal{F}$ of subspaces in $\F_q^n$ is called an $[n,k,L]_q$-\textit{weak subspace design} (c.f.~\cite{guruswami2016explicit}) if every $k$-dimensional subspace in $\F_q^n$ intersects non-trivially at most $L$ subspaces from $\mathcal{F}$. Despite it is not required by definition, many known constructions of weak subspace designs contain subspaces with a \textit{fixed} co-dimension at least $k$.  Weak subspace designs, almost sparse and almost affinely disjoint families have the following (trivial) relations:
\begin{itemize}
    \item an $[n,k,L]_q$-AS family is also an $[n,k+1,L]_q$-weak subspace design and an $[n,k,L-1]_q$-AAD family;
    \item for $n\ge 2k+1$, a partial $k$-spread of $\F_q^n$ is an $[n,k+1,L]_q$-weak subspace family if and only if it is also an $[n,k,L]_q$-AS family.
\end{itemize}
}
By explicit constructions presented in~\cite{guruswami2016explicit}, we derive that $p^{AS}(n,k,L)\ge \left \lfloor \frac{n-k}{k+1} \right \rfloor$ for $L\ge \frac{(n-1)(k+1)}{\lfloor(n-k)/(k+1)\rfloor}$. We note that the motivating application for subspace designs has $n$ growing and, thus, a straightforward application of these results does not give an optimal result for our problem except for the case $n=2k+1$. Constructions of subspace designs have found several  applications, such as constructing list-decodable error-correcting codes~\cite{guruswami2013list}, rank-metric codes~\cite{guruswami2013explicit}, dimensional expanders~\cite{guruswami2018lossless}. In particular, using this concept  the first deterministic polynomial time construction of list-decodable codes over constant-sized large alphabets and sub-logarithmic list size achieving the optimal rate has been designed~\cite{guruswami2016explicit}.   

 Almost affinely disjoint subspace families with $n=2k+1$ were first introduced in~\cite{polyanskii2019constructions}. It was proposed to use such families for constructing \textit{primitive batch codes} defined by Ishai et al.~in~\cite{ishai2004batch}. A primitive binary $[N,K,s]$-batch code encodes a binary string $\x$ of length $K$ into a binary string $\y$ of length $N$, such that each multiset of $k$ symbols
from $\x$ has $s$ mutually disjoint recovering sets from $\y$. The basic question is how to minimize the redundancy, $N-K$, for given parameters $K$ and $s$. Suppose that an $[n,k,L]_q$-AAD  family $\cF$ is given. Let $K:=q^n$. To construct a systematic batch code, we associate $K$ information bits with $K$ points in $\F_q^n$ and set $y_i:=x_i$ for $i\in[K]$. For an affine subspace of the form $\v+S$ with $S\in \cF$ and $\v\in\F_q^n$, we define a parity-check bit as a sum of information bits lying in this affine
subspace. As the number of distinct affine subspaces of such a form is $|\cF|q^{n-k}$, the constructed systematic code has  length $N=q^n+|\cF|q^{n-k}$ and the redundancy is $N-K=|\cF|q^{n-k}$. Moreover, it can be shown that this construction is an $[N,K,s]$-batch code with $s:=\lfloor |\cF|/L\rfloor$.

A naive way for constructing AAD families is by exploiting constructions of long linear codes with fixed \new{minimum} distance. Suppose that $H(\cC)$ is a parity-check matrix of a $q$-ary linear $[N,K,d]_q$-code $\cC$  of length $N$ and dimension $K$ with \new{minimum} distance $d=3k+1$. Let the subspace $S_i$ be the linear span of $k$ consecutive columns, from the $(ik+1)$-th to the $(i+1)k$-th column, of $H(\cC)$. Then $\cF:=\{S_1,\ldots,S_{\lfloor N/k\rfloor}\}$ is an $[N-K,k,1]_q$-AAD family. Thus, for fixed $N-K$, the longer the code, the larger the constructed family. Yekhanin and Dumer have developed a class of long non-binary codes with a fixed distance~\cite{yekhanin2004long}. For $k=1$ and $d=3k+1$, linear $[N,K,4]_q$-codes are known to be equivalent to caps in projective geometries and have been studied extensively under this name~\cite{hirschfeld2001packing,mukhopadhyay1978lower,edel1999recursive}. By the results of~\cite{yekhanin2004long,edel1999recursive}, for fixed $k$ and large enough $n$, it holds that $p^{ADD}(n,k,1)\ge (3k-1)(n+1)/(9k^2-9k+1)$.
\subsection{Our contribution}
The main results of our paper are presented in Theorems~\ref{th::upper bound AAD family}-\ref{th:general construction} showing that
$$
n-2k-\frac{(k+1)(n-k)}{L+1}\le p^{AS}(n,k,L)\le p^{AAD}(n,k,L-1)\le n-2k
$$
for all ranges of parameters. Moreover, by providing explicit constructions of AAD and AS subspace families we prove that $p^{AS}(n,1,L)=n-2$ for $L\ge n$ and $p^{AAD}(n,2,L)=n-4$ for $L\ge 4n^2-18n+21$, respectively.  Note that for the special case $n=2k+1$, by constructions from~\cite{guruswami2013explicit,polyanskii2019constructions} we have $p^{AAD}(2k+1,k,L)=1$ for $L\ge k$ and $p^{AS}(2k+1,k,L)=1$ for $L\ge 2k(k+1)$.

The remaining part of this paper is structured as follows. In Section~\ref{ss::converse bound}, we prove a non-existence result for AAD families. Random and explicit constructions for AS and AAD families
are presented in Section~\ref{ss::achievability}. We close our paper with a conjecture in
Section~\ref{ss::conclusion}.
\section{Converse bound}\label{ss::converse bound}
\begin{theorem}\label{th::upper bound AAD family}
Fix arbitrary positive integers $k$, $n$ and $L$ such that $2k< n$. Let $\cF$ be an $[n,k,L]_q$-AAD family. Then
\begin{equation}\label{eq::upper bound}
|\cF|\le 1+  L\frac{q^{n-k}-1}{q^{k}-1}.
\end{equation}
For $L=q^{o(1)}$, it follows that $p^{AS}(n,k,L)\le n- 2k$.
\end{theorem}
\begin{proof}
  By $m$ denote the cardinality of $\cF$ and let $S_i$ be the $i$th subspace from $\cF$.	 For  some $i\in[m]$, let $G_i\in \F^{n\times k}_q$ be a matrix whose columns span $S_i$ and let $H_i\in \F^{(n-k)\times n}_q$ be a matrix whose rows span $S_i^{\perp}$. Notice that $H_mG_j\in \F^{( n-k)\times k}_q$ has full column rank because $S_m$ and $S_j$ have \new{only} trivial intersection \new{for any $j\in [m-1]$, since they are from the AAD family which is a partial spread according to Definition~\ref{def::super nice def}.} Let $\hat G_j\in\Fq^{(n-2k)\times(n-k)}$ be a matrix whose rows form a basis  of the space orthogonal to the column span of $H_mG_j$.

  Next, we shall prove that for any non-zero vector $\w\in \F^{n-k}_q$, $\hat G_j\w$ is the all-zero vector for at most $L$  different $j$'s. To see this, suppose that for some multiset $\{j_1,\ldots, j_{L+1}\}\subset [m-1]$ of size $L+1$, we have $\hat G_{j_t} \w = \0$ \new{for every $t \in [L + 1]$}. This implies that $\w$ is in the column span of $H_mG_{j_t}$, say that $\w = H_mG_{j_t}\y_{t}$ for some $\y_{t}\in \F^ k_q$. Define $\v:=G_{j_1} \y_1$. So we have $$\w=H_mG_{j_t}\y_{t} = H_m\v$$ which means that $$G_{j_t}\y_t = \v + G_m\x_t$$ for some $\x_t\in\F_q^{k}$  (e.g., $\x_1 = \0$). But this implies  that $\v+S_m$ and $S_{j_t}$ intersect. By Definition~\ref{def::super nice def}, there are at most $L$ different $j$'s so that $\v+S_m$ and $S_j$ intersect. This leads to a contradiction.

  Whenever we have a collection of $m$ subspaces with the desired property, we have $(m-1)$ matrices $\hat G_j\in \F^{(n-2k)\times(n-k)}_q$ so that the span of any $L + 1$ of them has rank $n-k$. We claim that we must have \eqref{eq::upper bound} for such a collection of matrices $\{\hat G_j,\ j\in[m-1]\}$ to exist. Indeed, let $\x$ be a random non-zero vector in $\F^{n-k}_q$ . We observe that the expectation is $$\new{\mathbb{E}}\left|\{j\in[m-1]:\ \hat G_j \x = \0\}\right| =(m-1) \Pr \left\{\hat G_j\x = \0\right\} =(m-1)\frac{q^{k}-1}{q^{n-k}-1}.$$
  So if $m > L \frac{q^{n-k}-1}{q^k-1} + 1$, there exists some vector $\w \in \F_q^{n-k}$ so that $\hat G_j \w =\0$ for at least $L+1$ different $j$'s. This contradiction completes the proof.
\end{proof}
\begin{remark}
  Note that if we change the definition of an AAD subspace family by dropping the first property in Definition~\ref{def::super nice def}, then the matrices $H_m G_j$ would have full rank for at least $m-L-1$ different $j$'s. This results in the bound $m \le L \frac{q^{n-k}-1}{q^k-1} + L + 1.$
\end{remark}
\section{Constructions}\label{ss::achievability}
In this section, we provide new random and explicit constructions of AAD and AS families. We first show the existence result on AS families in Section~\ref{ss::random construction}. Section~\ref{ss::explicit construction} presents a novel construction \new{of AAD subspaces family based on Reed-Solomon codes} for $k=1$ and $k=2$.
\subsection{Random construction}
\label{ss::random construction}
	\begin{theorem}\label{th::probabillistic lower bound}
	For any fixed integers $L$, $n$, $k$ and $q\to\infty$, there exists an $[n,k,L]_q$-AS family of size $\underline{m}^{*}_q(n,k,L)$, where 
	$$
	\underline{m}^*_q(n,k,L):=q^{n-2k-\frac{(n-k)(k+1)}{(L+1)}}(1+o(1)).
	$$
For fixed $L$, it follows that $p^{AS}(n,k,L)\ge n-2k - (n-k)(k+1)/(L+1)$.
\end{theorem}
\begin{proof}
	The number of $k$-dimensional subspaces in $\F_q^n$ equals
	$$
	s_k:=\frac{(q^n-1)\dots(q^n-q^{k-1})}{(q^k-1)\dots (q^k-q^{k-1})} = \Theta(q^{k(n-k)}).
	$$
	We form a family of $k$-dimensional subspaces, written as $\cF:=\{S_1,\dots, S_M\}$, of size $M=q^{n-2k-(n-k)(k+1)/(L+1)}$ by choosing each subspace $S_i$ independently and equiprobably with probability $1/s_k$. So, it is possible that $S_i=S_j$ for some $i\neq j$.

	Define $\xi=|\{(i,j):\,i,j\in[M],\,i<j, |S_i\cap S_j|\neq 1\}|$, the number of pairs $(i,j)$ with $i< j$ such that $S_{i}$ and $S_{j}$ have the non-trivial intersection. We will estimate the expectation of $\xi$. The number of $k$-dimensional subspaces that do not intersect with a fixed  $k$-dimensional subspace (except the origin point) is equal to
	$$
	g_k:=\frac{(q^n-q^k)\dots(q^n-q^{2k-1})}{(q^k-1)\dots (q^k-q^{k-1})}.
	$$
	Thus, two random $k$-dimensional subspaces have a trivial intersection with probability $g_k/s_k$. The mathematical expectation of $\xi$ is then upper bounded as follows

	\begin{align*}
	\mathbb{E}(\xi) &\le \sum_{1\le i<j\le M}\Pr\{|V_i\cap V_j|\neq 1\} =  \sum_{1\le i<j\le M}\left(1-\frac{g_k}{s_k}\right)\\
	&\le\binom{M}{2}\left(1-\frac{(q^n-q^k)\dots(q^n-q^{2k-1})}{(q^n-1)\dots(q^n-q^{k-1})}\right)
	\\
	&<M^2\left.\left(1-\left(q^{nk}-q^{n(k-1)}\sum\limits_{i=k}^{2k-1}q^i\right)\right/q^{nk}\right)
	\\
	&=M^2(q^{2k-1-n}+o(q^{2k-1-n}))\new{<M(q^{-1}+o(q^{-1}))}.
	\end{align*}
	\new{By the Markov inequality
	$$
	\Pr(\xi>q^{0.5}\mathbb{E}(\xi))<\frac{\mathbb{E}(\xi)}{q^{0.5}\mathbb{E}(\xi)}=o(1).
	$$}
        \new{Since $q^{0.5}\mathbb{E}(\xi)<M(q^{-0.5}+o(q^{-0.5}))=o(M)$, we obtain that}
        with probability at least $1-o(1)$ there exists a family $\cF$ of size $M$, which contains at most $o(M)$ pairs of subspaces with the non-trivial intersection. If we delete one of the intersecting subspaces for each pair, then we obtain a family $\cF'\subset \cF$ of subspaces of size at least $M-o(M)$satisfying the first property of Definition~\ref{def::sparse family}.

	Now we compute the probability that the second property of Definition~\ref{def::sparse family} is violated. The number of $k$-dimensional subspaces that trivially intersect with a fixed  $(k+1)$-dimensional subspace, written as $V$, is equal to
	$$
	u_k:=\frac{(q^n-q^{k+1})\dots(q^n-q^{2k})}{(q^k-1)\dots (q^k-q^{k-1})}.
	$$
	Thus, the probability that $S_i$ does not intersect $V$ equals
	$u_k/s_k$. Let $\cF_V\subseteq\cF$ be the set of subspaces in $\cF$ that non-trivially intersect $V$, i.e., $\cF_V:=\{S\in\cF, |S\cap V|\neq 1\}$.  	Applying the union bound, we can estimate the probability that $V$ intersects at least $L+1$ subspaces in $\cF$ by
	\begin{align*}
	\Pr[|\cF_V|\geq L+1]&\leq \binom{M}{L+1}\left(1-\frac{u_k}{s_k}\right)^{L+1}\\
	&=\binom{M}{L+1}\left(1-\frac{(q^n-q^{k+1})\dots(q^n-q^{2k})}{(q^n-1)\dots(q^n-q^{k-1})}\right)^{L+1}\\
	&<M^{L+1}\left(1-\left(q^{nk}-q^{n(k-1)}\sum_{i=k+1}^{2k}q^{i}\right)/q^{nk}\right)^{L+1}
	\\
	&= M^{L+1}(q^{2k-n}+o(q^{2k-n}))^{L+1}\\
	&<q^{-(n-k)(k+1)}(1+o(1)).
	\end{align*}
	Recall that the total number of $(k+1)$-dimensional subspaces is $s_{k+1}$, which is $\Theta(q^{(k+1)(n-k-1)})$. Hence, by the union bound, $$s_{k+1}\cdot \Pr[|\cF_V\geq L+1|]<q^{-1-k}+o(q^{-1-k}),$$
	we have that with probability $o(1)$ the second property is violated. This completes the proof of the existence of a $[n,k,L]$-AS family of size $M-o(M)$.
\end{proof}
\subsection{Explicit constructions}\label{ss::explicit construction}
\begin{construction}\label{constr:: general k}
Let $q\ge nk$, $m=q^{n-2k}$ and $\gamma$ be a primitive element of $\F_q$. For $i\in[m]$, define $S_i$ to be the span of vectors $\{\v_{i,1},\ldots,
\v_{i,k}\}$ with
$$
\v_{i,j}:=\begin{pmatrix} \e_j & \Gamma_{j}(\c_{i}) &  h_j(\c_i)\end{pmatrix}\quad\text{for }j\in[k],
$$
where $\e_j$ is the unit vector of length $k$ having one in the $j$th position, $\c_i$ is a codeword of an $[n-k-1,n-2k,k]_q$ Reed-Solomon code having the following parity-check matrix
\begin{equation}\label{eq::parity-check matrix}
H_{RS}:=\begin{pmatrix}
1 & 1 & 1 & \cdots & 1\\
1 & \gamma & \gamma^2 & \cdots & \gamma^{n-k-2}\\
\vdots & \vdots & \vdots & \ddots & \vdots \\
1 & \gamma^{k-2} & \gamma^{2(k-2)} & \cdots & \gamma^{(n-k-2)(k-2)}
\end{pmatrix}
\end{equation}
and the map $\Gamma_j: \F_q^{n-k-1}\to \F_q^{n-k-1}$ is defined by
$$
\Gamma_j(\x):=\begin{pmatrix} \gamma^{j-1}x_1 &\gamma^{2(j-1)}x_2 &\gamma^{3(j-1)}x_3 & \cdots & \gamma^{(n-k-1)(j-1)}x_{n-k-1}\end{pmatrix},
$$
and the function $h_j(\x):=\sum_{p=1}^{n-k-1}x_{p}^{(j-1)(n-k-1)+p+1}$.
Then we set $\cF_{n,k}$ to be the collection of $S_i$'s.
\end{construction}
\begin{theorem}\label{th:general construction}
The family $\cF_{n,k}$ from Construction~\ref{constr:: general k} is a partial spread of $k$-dimensional subspaces in $\F_q^n$. Moreover, for $k=1$ and $k=2$, $\cF_{n,k}$ is $[n,k,L_{n,k}]_q$-AAD with $L_{n,1}=n-1$ and $L_{n,2}=1+2(n-2)(2n-5)$.
\end{theorem}
\begin{proof}
    The linear span of vectors $\{\v_{i,1},\ldots, \v_{i,k}\}$ defines a $k$-dimensional subspace in $\F_q^n$ as the restriction of $\v_{i,j}$ to the first $k$ coordinates is $\e_j$. Suppose that $S_i$ and $S_j$ have the non-trivial intersection. Thus, the rank of the system of vectors $\{\v_{i_1},\ldots, \v_{i,k},\v_{j,1},\ldots, \v_{j,k}\}$ is at most $2k-1$. This yields  that the rank of the system of vectors $$\left\{\Gamma_1(\c_i-\c_j), \Gamma_2(\c_i-\c_j),\ldots, \Gamma_{k}(\c_i-\c_j)\right\}$$ is not full. \new{Denote by $c_{i,j}$ the $j$th entry of $\c_i$.} Since $\c_i-\c_j$ is a non-zero codeword of the Reed-Solomon code with minimum distance $k$, there exist  $k$ coordinates $p_1,\ldots,p_{k}\in[n-k-1]$ such that $u_{t}:=c_{i,p_t}-c_{j,p_t}\neq 0$ for $t\in[k]$. Thus, after restricting each vector of the system onto coordinates $p_1,\ldots, p_k$, the non-full rank property is equivalent to that the determinant
    $$
    \det \begin{pmatrix}
    u_1 & u_2 & \cdots & u_k\\
    \gamma^{p_1}u_1 & \gamma^{p_2} u_2& \cdots & \gamma^{p_k} u_k \\
    \vdots & \vdots & \ddots & \vdots \\
    \gamma^{p_1(k-1)}u_1 & \gamma^{p_2(k-1)} u_2& \cdots & \gamma^{p_k(k-1)} u_k
    \end{pmatrix}
    = \prod_{t=1}^{k} u_t \prod_{1\le s<r\le k} (\gamma^{p_r}-\gamma^{p_s})
    $$
    is zero. However, since $\gamma$ is a primitive element of the field $\Fq$ with $q\geq nk$ and all $u_t$'s are non-zero, the determinant cannot be zero, which contradicts the assumption that $S_i$ and $S_j$ intersect non-trivially. Thus, the family $\cF_{n,k}$ is indeed a partial spread.

    Suppose that for some $\v\in\F_q^n\not\in S_i$, the linear span of $\v$ and $S_i$, written as $V$, intersects more than $L_{n,k}$ subspaces from the family $\cF_{n,k}$. Note that we can replace the vector $\v$ with a non-zero vector from one of the $L_{n,k}+1$ subspaces intersecting $V$. So, \new{we can assume that the vector $\v$ can be represented by $\v=\v_j(\aalpha):=\sum_{t=1}^k \alpha_t \v_{j,t}$ for some $\aalpha\in\F_q^{k}\setminus\{\0\}$, $j\in[m]\setminus \{i\}$}. In what follows, we estimate the number of $\ell\in[m]\setminus\{i,j\}$ such that there exists a $\bbeta\in\F_q^k$ and $\v_\ell(\bbeta)$ belongs to $V$. This is equivalent to the property that the system of vectors $\{\v_{i,1},\ldots, \v_{i,k}, \v_j(\aalpha), \v_\ell(\bbeta)\}$ is of rank at most $k+1$. By the structure of the vectors $\v_{i,t}$'s, this implies that the rank of
    $$
    R_{\bbeta,\ell}:= \begin{pmatrix}\sum_{t=1}^{k}\alpha_t \Gamma_t(\c_j-\c_i) & \sum_{t=1}^{k}\alpha_t (h_t(\c_j)-h_t(\c_i))\\
    \sum_{t=1}^{k}\beta_t \Gamma_t(\c_\ell-\c_i) & \sum_{t=1}^{k}\beta_t(h_t(\c_\ell)-h_t(\c_i))
    \end{pmatrix}
    $$
    is one.
    \new{Observe that the $p$th element, $p\in[n-k-1]$, of the first row of $R_{\bbeta,\ell}$ has the form $(c_{i,p}-c_{j,p})\sum_{t=1}^{k}\alpha_t\gamma^{(t-1)p}$.
    We can think about $\sum_{t=1}^{k}\alpha_t\gamma^{(t-1)p}$ as of a polynomial of degree $k-1$ in a variable $x=\gamma^p$. It has at most $k-1$ roots in $\F_q$. For each root $x_0$, there exists at most one $p\in[n-k-1]$ such that $\gamma^p=x_0$, since $\gamma$ is a primitive element in $\Fq$ with $q\geq nk$. Therefore, for any nonzero $\aalpha$, there are at most $k-1$ different $p\in[n-k-1]$ so that $\sum_{t=1}^{k}\alpha_t\gamma^{(t-1)p}=0$. Since $\c_i$ and $\c_j$ are codewords of a linear code with minimum distance $k$, there are at least $k$ positions $p\in[n-k-1]$ such that $c_{i,p}-c_{j,p}\neq 0$. Hence, there is at least one position $p_0\in[n-k-1]$ such that the $p_0$th entry of the vector $\sum_{t=1}^{k}\alpha_t \Gamma_t(\c_j-\c_i)$ is non-zero.
    }
    \begin{lemma}\label{lem::given beta number of solutions}
    Given a non-zero vector $\bbeta\in\F_q^k$, there exist at most $k(n-k)$ different $\c_\ell$'s from the Reed-Solomon code defined by~\eqref{eq::parity-check matrix} such that the matrix $R_{\bbeta,\ell}$ has rank $1$.
    \end{lemma}
    \begin{proof}
        If $R_{\bbeta,\ell}$ has rank $1$, then every column in the left part of $R_{\bbeta,\ell}$ (consisting of the first $n-k-1$ columns) is linearly dependent of the $p_0$th column, which is known to be non-zero. Moreover, the dependence can be found by the first row in $R_{\bbeta,\ell}$. Thus for any $p\in[n-k-1]$, there exists some element $\phi_p\in\F_q$ so that
        \begin{equation}\label{eq::first system}
             (c_{\ell,p} - c_{i,p})\sum_{t=1}^{k}\beta_t \gamma^{(t-1)p} =  \phi_p(c_{\ell,p_0} - c_{i,p_0})\sum_{t=1}^{k}\beta_t \gamma^{(t-1)p_0}.
        \end{equation}
        Recall that $\sum_{t=1}^{k}\beta_t \gamma^{(t-1)p}=0$ for at most $k-1$ values of $p\in[n-k-1]$.
        We have additionally restrictions imposed by the parity-check matrix~\eqref{eq::parity-check matrix}, namely, $\forall t\in[k-1]$
        \begin{equation}\label{eq::second system}
        \sum_{p=1}^{n-k-1} \gamma^{(p-1)(t-1)}(c_{\ell,p} - c_{i,p}) = 0.
        \end{equation}
        Thus, the system of equations~\eqref{eq::first system}-\eqref{eq::second system} with at least $n-k-2$ linear independent equations for variables $\{c_{\ell,p},p\in[n-k-1]\setminus\{p_0\}\}$ has at most one solution. In other words,  $c_{\ell,p}=a_p c_{\ell,p_0}+b_p$ with some $a_p,b_p\in\F_q$ and $p\in[n-k-1]$. This implies that the determinant
        $$
        \det
        \begin{pmatrix}\sum_{t=1}^k \alpha_t \gamma^{(t-1)p_0}(c_{j,p_0}-c_{i,p_0}) & \sum_{t=1}^{k}\alpha_t(h(\c_j)-h(\c_i))\\
        \sum_{t=1}^k \beta_t \gamma^{(t-1)p_0}(c_{\ell,p_0}-c_{i,p_0}) & \sum_{t=1}^{k}\beta_t(h(\c_\ell)-h(\c_i))
        \end{pmatrix}
        $$
        is zero. Note that the entry
        \begin{align*}&\sum_{t=1}^{k}\beta_t(h(\c_\ell)-h(\c_i)) \\
        &= \sum_{t=1}^{k}\beta_t \sum_{p=1}^{n-k-1}\left((a_p c_{\ell,p_0}+b_p)^{(t-1)(n-k-1)+p+1}-c_{i,p}^{(t-1)(n-k-1)+p+1}\right)
        \end{align*}
        is a polynomial in $c_{\ell,p_0}$ of degree at least $p_0+1$ and at most $k(n-k-1)+1\le k(n-k)$. Therefore, the determinant represents a non-trivial univariate polynomial in $c_{\ell,p_0}$ of degree at least $p_0+1$ and at most $k(n-k)$. Since $q\ge nk$, there are at most $k(n-k)$ solutions for $c_{\ell,p_0}$ vanishing the determinant.
    \end{proof}
    Let us proceed with proving the remaining statement of this theorem.  For this purpose, for $k=1,2$, we estimate the number of possible $\bbeta$'s such that the first $n-k-1$ columns of $R_{\beta,\ell}$, treated as vectors, are \new{parallel} to the $p_0$th column and apply Lemma~\ref{lem::given beta number of solutions}.

    For the case $k=1$, we first note that we can scale vector $\bbeta$ and, thus, have $\bbeta = (\beta_1)=(1)$. By Lemma~\ref{lem::given beta number of solutions}, the number of distinct appropriate $\ell$'s is at most $n-1$. Thus, $\cF_{n,1}$ is an $[n,1,L]_q$-AAD family.

    Now we discuss the case $k=2$. Suppose that $\bbeta$ is not \new{parallel} to the vector $(1,-1)$.
    Since $\bbeta$ can be appropriately scaled, we shall estimate the number of \textit{suspicious} $\bbeta=(\beta_1,1-\beta_1)$ which means that there might exist some $\c_\ell$ from the Reed-Solomon code for such $\bbeta$ so that $R_{\bbeta,\ell}$ is degenerate. Define the set $B$ that includes $\bbeta = (1,0)$ and $\bbeta$'s with the property $\beta_1+(1-\beta_1) \gamma^{p}=0$ for some $p\in[n-3]$. Consider a $\bbeta\not\in B$. If $R_{\bbeta,\ell}$ has rank $1$, then two rows of  $R_{\bbeta,\ell}$, treated as vectors, are \new{parallel} and there exists some non-zero $\lambda\in\F_q$ such that
    \begin{align}\label{eq:colinear-constraint}
        c_{\ell,p}-c_{i,p}=\lambda\ \frac{(\alpha_1+\alpha_2\gamma^{p})(c_{j,p}-c_{i,p})}{\beta_1+(1-\beta_1)\gamma^p}\quad\text{for }p\in[n-3].
    \end{align}

    Let $w_p$ denote the numerator of the above fraction. Recall that there exists  $p_0$ such that $w_{p_0}\neq 0$. From the parity-check property $\sum_{p=1}^{n-3}(c_{\ell,p}-c_{i,p})=0$ imposed by~\eqref{eq::parity-check matrix}, we thus derive
    $$
    \sum_{p=1}^{n-3} \frac{w_p}{\beta_1+(1-\beta_1)\gamma^p} = 0\quad \Leftrightarrow  \quad \sum_{p=1}^{n-3}w_p \prod_{t\neq p}(\beta_1+(1-\beta_1)\gamma^t)=0.
    $$
     We can think of the left-hand side of the above equation as a non-trivial univariate polynomial in $\beta_1$ of degree at most $n-4$. Indeed, the polynomial is non-trivial as its evaluation at $\beta_1^*$ satisfying $\beta_1^* +(1-\beta_1^*)\gamma^{p_0}=0$ is $w_{p_0}\prod\limits^{n-3}_{\substack{t=1\\t\neq p_0}}(\beta_1^*+(1-\beta_1^*)\gamma^t)\neq 0$.
    Therefore, there are at most $n-4$ suspicious $\bbeta$'s not included to $B$ and not \new{parallel} to $(1,-1)$ such that the vector $\sum_{t=1}^{2}\beta_t \Gamma_t(\c_j-\c_i)$ is \new{parallel} to $\sum_{t=1}^{2}\beta_t \Gamma_t(\c_\ell-\c_i)$. Define $D$  to be the the union of suspicious $\bbeta$'s, the set $B$ and the vector $(1,-1)$. As $|D|\le 2n-5$, it can be easily verified that
    by Lemma~\ref{lem::given beta number of solutions}, $\cF_{n,2}$ is an $[n,2,L]_q$-AAD family with $L=1+2(n-2)(2n-5)$.
    \end{proof}
\begin{remark}
\new{Note that any $2$-dimension subspace, i.e., a plane, that intersects with a $1$-dimensional subspace, i.e., a line, from an $[n,k=1,L]_q$-ADD family, must contain the line. Therefore,} any $[n,1,L]_q$-AAD family is also an $[n,1,L+1]_q$-AS family.
Therefore, $p^{AS}(n,1,n)=p^{\new{AAD}}(n,1,n-1)=n-2$.
\end{remark}

\section{Conclusion}\label{ss::conclusion}
In this paper, we have introduced new notions concerning vector finite spaces which we have called almost affinely disjoint and almost sparse families of subspaces. We have presented lower and upper bounds on the polynomial growth of the maximal sizes of these families. For the cases $k=1$ and $k=2$, our explicit constructions asymptotically achieve the converse bound. We conjecture that for any $k,n$ and a large enough $L=L(n,k)$, the quantities $p^{ADD}(n,k,L)$ and $p^{AS}(n,k,L)$ are equal to $n-2k$.

\section{Acknowledgement}
H. Liu was funded by the German Israeli Project Cooperation (DIP) grant under grant no.~KR3517/9-1.
N. Polyanskii was funded by the German Research Foundation (Deutsche Forschungsgemeinschaft, DFG) under Grant No. WA3907/1-1.
Ilya Vorobyev was funded by RFBR and JSPS, project number 20-51-50007, and by RFBR, project number 20-01-00559.
\bibliographystyle{elsarticle-harv}
\bibliography{refs}
\end{document}